\documentclass[11pt]{article}
\usepackage[dvips]{graphics}
\usepackage{multirow}
\usepackage{array}
\usepackage{epsfig,rotating}
\usepackage{amssymb,amsmath}
\usepackage{latexsym}
\usepackage[usenames]{color}

\numberwithin{equation}{section}

\newcommand{\be}{\begin{equation}}
\newcommand{\ee}{\end{equation}}
\newcommand{\beaa}{\begin{eqnarray*}}
\newcommand{\eeaa}{\end{eqnarray*}}
\newcommand{\bea}{\begin{eqnarray}}
\newcommand{\eea}{\end{eqnarray}}

\newcommand{\bei}{\begin{itemize}}
\newcommand{\eei}{\end{itemize}}

\newcommand{\bd}{\bold}

\newtheorem{theorem}{ \noindent T{\footnotesize HEOREM}}

\newtheorem{lemma}{ \noindent L{\footnotesize EMMA}}[section]

\newcommand\td{\overset{d}{\to}}
\begin{document}
%\nocite{*}

\title{Spectral Radii of Truncated Circular Unitary Matrices}
\author{Wenhao Gui$^{1}$
and  Yongcheng Qi$^2$\\
Beijing Jiaotong University and University of Minnesota Duluth}

\date{}
\maketitle

\footnotetext[1]{Department of Mathematics, Beijing Jiaotong University, Beijing 100044, China %\newline  \indent \  \
}

\footnotetext[2]{Department of Mathematics and Statistics, University of Minnesota Duluth, MN 55812, USA, yqi@d.umn.edu.
%\newline  \indent \ \
%The research of Yongcheng Qi was supported in part by NSF Grant DMS-1005345.
}

\begin{abstract}
\noindent  Consider a truncated circular unitary matrix which is a $p_n$ by $p_n$ submatrix of an $n$ by $n$ circular
unitary matrix by deleting the last $n-p_n$ columns and rows.   Jiang and Qi (2017) proved that the maximum absolute value of the eigenvalues (known as spectral radius) of the truncated matrix, after properly normalized, converges in distribution to the Gumbel  distribution if $p_n/n$ is bounded away from $0$ and $1$.  In this paper we investigate the limiting distribution of the spectral radius under one of the following four conditions: (1). $p_n\to\infty$ and $p_n/n\to 0$ as $n\to\infty$; (2). $(n-p_n)/n\to 0$ and $(n-p_n)/(\log n)^3\to\infty$ as
$n\to\infty$; (3). $n-p_n\to\infty$ and $(n-p_n)/\log n\to 0$ as $n\to\infty$ and (4). $n-p_n=k\ge 1$ is a fixed integer.  We prove that the spectral radius converges in distribution to the Gumbel distribution under the first three conditions and to a reversed Weibull distribution  under the fourth condition.
\end{abstract}

%{\red [see, e.g., Theorem 2.8.1
% on page 110 from Andrews, Askey and Roy (1999)]}\\\\


%\noindent\textbf{AMS 2000 Subject Classification: \/}  15A52, 60F99, 60G55, 60G70. \\


\newpage

\section{Introduction}\label{intro}

The early study of large random matrices was stimulated by analysis of high-dimensional data. One example is Wishart's (1928)
investigation on large covariance matrices whose statistical properties are mainly determined by
eigenvalues and eigenvectors from the point view of a principal components analysis.   Since then, the random matrix theory has been
developed very rapidly and found many applications in areas such as heavy-nuclei atoms (Wigner, 1955),
number theory (Mezzadri and Snaith, 2005), quantum mechanics (Mehta, 2005), condensed matter physics (Forrester, 2010),
wireless communications (Couillet and Debbah, 2011).

The study of random matrices has greatly been motivated by Tracy and Widom's (1994, 1996) work. They show that
the largest eigenvalues of the three Hermitian matrices (Gaussian orthogonal ensemble, Gaussian unitary ensemble and Gaussian symplectic ensemble) converge to some special distributions that are now known as the Tracy-Widom laws.
Subsequently, the Tracy-Widom laws have found their applications in the study of problems such as
the longest increasing subsequence (Baik et al., 1999), combinatorics, growth processes, random tilings and the determinantal point processes  (see, e.g., Tracy and Widom (2002), Johansson (2007) and references therein) and the largest eigenvalues in the high-dimensional statistics
(see, e.g.,  Johnstone (2001, 2008) and Jiang (2009)).
 Some recent research focuses on  the universality of the largest eigenvalues of matrices with non-Gaussian entries; see, for example, Tao and Vu (2011), Erd\H{o}s et al. (2012) and the references therein.

 Consider a non-Hermitian matrix $\bd{M}$ with eigenvalues $z_1, \cdots, z_n$.  The largest absolute values of the eigenvalues
 $\max_{1\leq j \leq n}|z_j|$ is refereed to as the spectral radius of $\bd{M}$. The spectral radii of the real, complex and symplectic Ginibre ensembles are investigated by Rider (2003, 2004) and Rider and Sinclair (2014), and  it is proved that the spectral radius for the complex Ginibre ensemble converges to the Gumbel distribution.  This indicates that non-Hermitian matrices exhibit quite different
 behaviors from Hermitian matrices in terms of the limiting distribution for the largest absolute values of the eigenvalues.

 A very recent paper by Jiang and Qi (2017) studies the largest radii of three rotation-invariant and non-Hermitian random matrices: the spherical ensemble, the truncation of circular unitary ensemble and the product of independent complex Ginibre ensembles. It is proved in the paper that the spectral radii converge to the Gumbel distribution and some new distributions.

 The circular unitary ensemble is an $n\times n$ random matrix with Haar measure on the unitary group, and it is also called Haar-invariant unitary matrix. Let $\bd{U}$ be an $n\times n$ circular unitary matrix. The $n$ eigenvalues of the circular unitary matrix $\bd{U}$ are distributed over $\{z\in\mathcal{C}: |z|=1\}$ , where $\mathcal{C}$ is the complex plane, and their joint density function is given by
\[
\frac{1}{n!(2\pi)^n}\cdot\prod_{1\leq j<k\leq n}|z_j - z_k|^2;
\]
see, e.g., Hiai and Petz (2000).

For $n>p\geq 1$, write
\begin{eqnarray*}
\bd{U}=\begin{pmatrix}
\bd{A}\ \ \bd{C}^*\\
\bd{B}\ \ \bd{D}
\end{pmatrix}
\end{eqnarray*}
where $\bd{A}$, as a truncation of $\bd{U}$, is a $p\times p$ submatrix. Let $z_1, \cdots, z_p$ be the eigenvalues of $\bd{A}$. Then their density function is
\begin{equation}\label{Good_forever}
C\cdot\prod_{1\leq j<k\leq p}|z_j - z_k|^2\prod_{j=1}^{p}(1-|z_j|^2)^{n-p-1}
\end{equation}
where $C$ is a normalizing constant. See, e.g., Zyczkowski and Sommers (2000).

Assume $p=p_n$ depends on $n$ and set $c=\lim_{n\to\infty}\frac{p_n}{n}$.  \.{Z}yczkowski and Sommers (2000) show that the empirical distribution of $z_i$'s converges to the distribution with density proportional to $\frac{1}{(1-|z|^2)^2}$ for $|z|\leq c$ if $c\in (0,1)$.
Dong et al. (2012) prove that the empirical distribution goes to the circular law and the arc law as $c=0$ and $c=1$, respectively.
See also Diaconis and Evans (2001) and Jiang (2009, 2010) and references therein for more results.

Jiang and Qi (2017) have proved that the spectral radius $\max_{1\le j\le p}|z_j|$ for the truncated circular unitary ensemble
converges to the Gumbel distribution when the dimension of the truncated
truncated circular unitary matrix is of the same order as the dimension of the original circular unitary matrix, see Theorem~\ref{thm1} in section~\ref{main}.

In this paper we consider heavily truncated and lightly truncated circular unitary matrices and investigate
the limiting distribution of the spectral radii for those truncated circular unitary matrices.
Our results complement that in Jiang and Qi (2017).

The rest of the paper is organized as follows.  The main results in this paper are given in section~\ref{main} and their proofs
 are provided in section~\ref{proof}.

%including three theorems on the limiting distributions of the spectral radii of truncated unitary ensembles corresponding to

\section{Main Results}\label{main}

Consider the $p_n\times p_n$ submatrix $\mathbf{A}$, truncated from a $n\times n$ circular unitary matrix $\mathbf{U}$
in section~\ref{intro}.
Denote the $p_n$ eigenvalues  as $z_1,\cdots, z_{p_n}$ with the joint density function given by \eqref{Good_forever}.

For completeness, we first quote a theorem in Jiang and Qi (2017) on the limiting distribution of the spectral radii $\max_{1\le j\le p_n}|z_j|$ before we give our results in the paper.

\begin{theorem}\label{thm1}
Assume that $z_1, \cdots, z_{p_n}$ have density as in (\ref{Good_forever}) and there exist constants $h_1, h_2\in (0,1)$ such that $h_1< \frac{p_n}{n}< h_2$ for all $n\geq 2.$ Then $(\max_{1\leq j \leq p_n}|z_j|-A_n)/B_n$ converges weakly to the Gumbel distribution $\Lambda(x)=\exp(-e^{-x})$, $x\in \mathbb{R}$,
where $A_n=c_n+\frac{1}{2}(1-c_n^2)^{1/2}(n-1)^{-1/2}a_n$,  $B_n=\frac{1}{2}(1-c_n^2)^{1/2}(n-1)^{-1/2}b_n$,
\[
c_n=\Big(\frac{p_n-1}{n-1}\Big)^{1/2},~~ b_n=b\Big(\frac{nc_n^2}{1-c_n^2}\Big), ~~a_n=a\Big(\frac{nc_n^2}{1-c_n^2}\Big)
\]
with
%\begin{equation}\label{abx}
\beaa
a(y)=(\log y)^{1/2}-(\log y)^{-1/2}\log(\sqrt{2\pi}\log y)\ \mbox{ and } \ b(y)=(\log y)^{-1/2}
\eeaa
%\end{equation}
for $y>3$.
\end{theorem}

Note that in Theorem~\ref{thm1}, a restriction on the dimension $p_n$ of the truncated circular unitary matrix
$\bd{A}$ is made as follows: there exist some $0<h_1<h_2<1$ such that $h_1n<p_n<h_2n$ for all large $n$.
In this paper we are devoted to study of the spectral radii
$\max_{1\le j\le p_n}|z_j|$ in the following conditions:

\begin{equation}\label{c1}
p_n\to\infty \mbox{ and } \frac{p_n}{n}\to 0\mbox{ as }n\to\infty;
\end{equation}

\begin{equation}\label{c2}
\frac{n-p_n}{(\log n)^3}\to\infty\mbox{ and } \frac{n-p_n}{n}\to 0\mbox{ as }n\to\infty;
\end{equation}

\begin{equation}\label{c3}
n-p_n\to\infty \mbox{ and } \frac{n-p_n}{\log n}\to 0\mbox{ as }n\to\infty;
\end{equation}

\begin{equation}\label{c4}
n-p_n=k\ge 1 \mbox{ is fixed integer.}
\end{equation}

%%%where $k_n=n-p_n$.

%%%Gumbel-type, Fr\'echet-type, and

The main results of the paper are the following theorems:

\begin{theorem}\label{thm2}
Under condition \eqref{c1} or \eqref{c2}, $(\max_{1\leq j \leq p_n}|z_j|-A_n)/B_n$ converges weakly to the Gumbel distribution $\Lambda(x)=\exp(-e^{-x})$, $x\in \mathbb{R}$, where $A_n$ and $B_n$ are defined as in Theorem~\ref{thm1}.
\end{theorem}

\begin{theorem}\label{thm3}
Under condition \eqref{c3},  $(\max_{1\leq j \leq p_n}|z_j|-A_n)/B_n$ converges weakly to the Gumbel distribution
$\Lambda(x)=\exp(-e^{-x})$, $x\in \mathbb{R}$, where $A_n=(1-a_n/n)^{1/2}$ and $B_n=a_n/(2nk_n)$, where $a_n$ is given by
\[
\frac{1}{(k_n-1)!}\int^{a_n}_0t^{k_n-1}e^{-t}dt=\frac{k_n}{n}.
\]
where $k_n=n-p_n$.
\end{theorem}

\begin{theorem}\label{thm4}
Under condition \eqref{c4},   $\frac{2n^{1+1/k}}{((k+1)!)^{1/k}}(\max_{1\leq j \leq p_n}|z_j|-1)$ converges weakly to
the reversed Weibull distribution $W_{k}(x)$ defined as
\[
W_{k}(x)=\left\{
                \begin{array}{ll}
                 \exp(-(-x)^{k}) , & \hbox{$x\le 0$;} \\
                  1, & \hbox{$x>0$.}
                \end{array}
              \right.
\]
\end{theorem}

We notice that the limiting distribution of the spectral radii depends on the dimension of truncated matrices.
Our results in Theorems~\ref{thm2}, \ref{thm3} and \ref{thm4} indicate that the limiting distribution of the spectral radii of the truncated
circular unitary matrices is Gumbel distribution $\Lambda$ if the parameter $k_n=n-p_n$, the number of truncated columns and rows diverges.
When the truncation is very light, that is, $k_n=n-p_n=k\ge 1$ is a fixed integer, the limiting distribution of the spectral radii of the truncated matrices is the reversed Weibull distribution $W_k$.

It is obvious that the case when $k_n=n-p_n$ is of order between $\log n$ and $(\log n)^3$ has not been covered in Theorems~\ref{thm1} to \ref{thm4}. We conjecture that $\max_{1\leq j \leq p_n}|z_j|$, after properly normalized, converges in distribution to the Gumbel distribution in this case.

\section{Proofs}\label{proof}

%In the rest of the paper, we will need the following notation. The symbol
%$C_n\sim D_n$ as $n\to\infty$ means that $\lim_{n\to\infty}\frac{C_n}{D_n}=1$. Similarly,  $C_n(t)\sim D_n(t)$ uniformly over $t\in
%T_n$ if $\lim_{n\to\infty}\sup_{t\in T_n }|\frac{C_n(t)}{D_n(t)}-1|=0$. Also,  $C_n(t)=O(D_n(t))$ uniformly over $t\in T_n$
%if $\sup_{t\in T_n}|\frac{C_n(t)}{D_n(t)}|$ is bounded. We write  $C_n(t)=o(D_n(t))$ uniformly over $t\in T_n$
%if $\sup_{t\in T_n}|\frac{C_n(t)}{D_n(t)}|$ converges to zero as $n\to\infty$.

Set $k_n=n-p_n$ and define $a(x)$ and  $b(x)$ as in Theorem \ref{thm1}.
Define $\phi(x)=\frac{1}{\sqrt{2\pi}}e^{-x^2/2}$  and $\Phi(x)=\frac{1}{\sqrt{2\pi}}\int_{-\infty}^xe^{-t^2/2}\,dt$ for $x\in \mathbb{R}$, the density function and the cumulative distribution of the standard normal, respectively. The symbol
$C_n\sim D_n$ as $n\to\infty$ means that $\lim_{n\to\infty}\frac{C_n}{D_n}=1$.

For random variables $\{X_n;\, n\geq 1\}$ and constants $\{a_n;\, n\geq 1\}$, we write $X_n=O_p(a_n)$ if $\lim_{x\to +\infty}\limsup_{n\to\infty}P(|\frac{X_n}{a_n}|\geq x)=0$. It is well known that $\frac{X_n}{a_nb_n}\to 0$ in probability as $n\to\infty$
if $X_n=O_p(a_n)$ and $\{b_n;\, n\geq 1\}$ is a sequence of constants with $\lim_{n\to\infty}b_n=\infty$.

Let $U_i$, $i\ge 1$ be a sequence of i.i.d. random variables uniformly distributed over $(0,1)$, and $U_{1:n}\le U_{2:n}\le \cdots\le U_{n:n}$ be the order statistics of $U_1, U_2, \cdots, U_n$ for each $n\ge 1$.
Then
from page 14 on the book by Balakrishnan and Cohen (1991), we know that the cumulative distribution function of $U_{i:n}$ is given by
\begin{equation}\label{cdf}
F_{i:n}(x)=\sum^n_{r=i}{{n}\choose{r}}x^r(1-x)^{n-r}=\frac{n!}{(i-1)!(n-i)!}\int^{x}_0
t^{i-1}(1-t)^{n-i}dt, ~0\le x\le 1
\end{equation}
 for each $1\le i\le n$, and the probability density function (pdf) of  $U_{i:n}$ is given by
  \begin{equation}\label{pdf}
f_{i:n}(x)=\frac{n!}{(i-1)!(n-i)!}
x^{i-1}(1-x)^{n-i}, ~~~0<x<1.
\end{equation}
This is the so-called Beta distribution, denoted by Beta($i, n-i+1$).

From \eqref{pdf},  $U_{p_n-j+1: n-j}$ has a Beta($p_n-j+1, k_n$) distribution with pdf given by
\[
f_{p_n-j+1:n-j}(x)=\frac{(n-j)!}{(p_n-j)!(k_n-1)!}x^{p_n-j}(1-x)^{k_n-1},~~~~~~x\in (0,1).
\]

For each $n\geq 2$, let $\{Y_{nj};\, 1\leq j \leq p_n\}$ be independent random variables such that
$Y_{nj}$ and $(U_{p_n-j+1: n-j})^{1/2}$ have the same distribution.  Jiang and Qi (2017) have shown that
$\max_{1\leq j \leq p_n}|z_j|$ and $\max_{1\leq j \leq p_n}Y_{nj}$ have the same distribution, that is
\begin{equation}\label{identical}
P\Big(\max_{1\leq j \leq p_n}|z_j|^2\le t\Big)=P\Big(\max_{1\leq j \leq p_n}Y_{nj}^2\le t\Big)=\prod^{p_n}_{j=1}F_{p_n-j+1:n-j}(t).
\end{equation}
for any $0<t<1$, and
\begin{equation}\label{monotone2}
1-F_{1:k_n}(x)\le 1-F_{2: k_n+1}(x)\le \cdots\le 1-F_{p_n: n-1}(x)
\end{equation}
for $x\in (0,1)$.  See the proof of Theorem 2 in Jiang and Qi (2017).

%%Beta and F distribution:
%%Beyer, W. H. CRC Standard Mathematical Tables, 28th ed. Boca Raton, FL: CRC Press, p. 536, 1987.

%%We will need some properties of Gamma distributions.  Assume $\alpha>0$ and $\beta>0$. The density of a Gamma($\alpha, \beta$)
%%distribution is given by
%%\[
%%g(x, \alpha, \beta)=\frac{1}{\beta^\alpha\Gamma(\alpha)}x^{\alpha-1}e^{-x/\beta},~~~x>0,
%%\]
%%where $\Gamma(\alpha)$ is the gamma function defined as
%%\[
%%\Gamma(\alpha)=\int^\infty_0x^{\alpha-1}e^{-x}dx.
%%\]

\subsection{Preliminary Lemmas}

We will present some useful lemmas before we prove our main results.

\begin{lemma}\label{prod2sum} Suppose $\{l_n;\, n\geq 1\}$ is sequence of positive integers. Let $z_{nj}\in [0,1)$ be real numbers for $1\leq j \leq l_n$ such that $\max_{1\le j\le l_n}z_{nj}\to 0$ as $n\to\infty$. Then $\displaystyle\lim_{n\to\infty}\prod^{l_n}_{j=1}(1-z_{nj})\in (0,1)$ exists if and only if the limit $\displaystyle\lim_{n\to\infty}\sum^{l_n}_{j=1}z_{nj}=:z\in (0,\infty)$ exists and the relationship of the two limits is given by
\begin{equation}\label{prod-sum}
\lim_{n\to\infty}\prod^{l_n}_{i=1}(1-z_{ni})=e^{-z}.
\end{equation}
\end{lemma}

\noindent{\bf Proof.}  From the Taylor expansion
\[
\log(1-x)=-x+O(x^2)~~~~\mbox{ as } x\to 0,
\]
which implies $\log(1-z_{nj})=-z_{nj}+O(z_{nj}^2)$ uniformly over $1\le j\le l_n$ since $\max_{1\le j\le l_n}z_{nj}\to 0$
as $n\to\infty$. Therefore,
\[
\prod^{l_n}_{j=1}(1-z_{nj})=\exp(\sum^{l_n}_{j=1}\log (1-z_{nj}))=\exp(-(1+O(\max_{1\le j\le l_n}z_{nj}))\sum^{l_n}_{j=1}z_{nj})
\]
The lemma can be easily concluded from the above expression. \hfill$\blacksquare$

\begin{lemma}\label{little} Let $\{l_n\}$ be a sequence of positive integers such that $l_n\to\infty$ and for each $n$, $\{z_{nj}, ~1\le j\le l_n\}$ are non-negative numbers such that $z_{nj}$ is non-increasing in $j$ with $z_{n1}>0$.  Then for any sequence of positive integers $\{r_n\}$satisfying that $r_n<l_n$ for all large $n$ and $r_n/l_n\to 1$ as $n\to\infty$, we have
\[
\frac{\sum^{l_n}_{j=1}z_{nj}}{\sum^{r_n}_{j=1}z_{nj}}\to 1
\]
as $n\to\infty$.
\end{lemma}

\noindent{\bf Proof.}  It suffices to show that
\begin{equation}\label{short}
\frac{\sum^{l_n}_{j=r_n+1}z_{nj}}{\sum^{r_n}_{j=1}z_{nj}}\to 0
\end{equation}
as $n\to\infty$.  In fact,  from the monotonicity of $z_{nj}$, we have $z_{nj}\le \frac{1}{r_n}\sum^{r_n}_{j=1}z_{nj}$ for $r_n+1\le j\le l_n$. Hence
\[
\sum^{l_n}_{j=r_n+1}z_{nj}\le \frac{l_n-r_n}{r_n}\sum^{r_n}_{j=1}z_{nj}
\]
which implies
\[
\frac{\sum^{l_n}_{j=r_n+1}z_{nj}}{\sum^{r_n}_{j=1}z_{nj}}\le \frac{l_n-r_n}{r_n}\to 0,
\]
proving \eqref{short}.    \hfill$\blacksquare$

For the rest of  the proofs,  define
\begin{equation}\label{znj}
z_{nj}=1-F_{p_n-j+1:n-j}(t_n), ~~~~1\le j\le p_n,
\end{equation}
where $t_n\in (0,1)$ will be specified later in the proof of each theorem. From \eqref{monotone2},
\begin{equation}\label{monotonez}
z_{n1}\ge z_{n2}\ge \cdots\ge z_{np_n}\ge 0.
\end{equation}
Obviously,  we have for $1\le j\le p_n$.
\bea\label{znj=}
z_{nj}&=&\frac{(n-j)!}{(p_n-j)!(k_n-1)!}\int^1_{t_n}t^{p_n-j}(1-t)^{k_n-1}dt\nonumber\\
&=&\frac{(n-j)!}{(p_n-j)!(k_n-1)!}\int^{1-t_n}_0(1-t)^{p_n-j}t^{k_n-1}dt.
\eea

\begin{lemma}\label{less} Assume that $1\le p_n<n$ and $p_n\to \infty$ as $n\to\infty$.
Let $\{r_n\}$ satisfy the condition in Lemma~\ref{little} with $l_n=p_n$.   Assume $\alpha_n>0$ and $\beta_n$ are real numbers such that
$\lim_{n\to\infty}P(Y_{n1}^2>\beta_n+\alpha_nx)= 0$ for any $x\in R$. If $(\max_{1\le j\le r_n}Y_{nj}^2-\beta_n)/\alpha_n$ converges in distribution to a
cdf $G$, then $(\max_{1\le j\le p_n}Y_{nj}^2-\beta_n)/\alpha_n$ converges in distribution to  the same distribution $G$.
\end{lemma}

\noindent{\bf Proof.} Note that $(\max_{1\le j\le r_n}Y_{nj}^2-\beta_n)/\alpha_n$ converges in distribution to the
cdf $G$ if and only if
\begin{equation}\label{cdfG-rn}
\lim_{n\to\infty}P(\max_{1\le j\le r_n}Y_{nj}^2\le \beta_n+\alpha_nx)=G(x)
\end{equation}
for every continuity point $x$ of $G$ with $G(x)\in (0,1)$. We need to prove the above expression is still true
when $r_n$ is replaced by $p_n$.  Now fix $x$, a continuity point of $G$ with $G(x)\in (0,1)$. Set $t_n=t_n(x)=\beta_n+\alpha_nx$ and
 define $z_{nj}$ as in \eqref{znj}. Note that \eqref{monotonez} holds,  $z_{n1}\to 0$ as $n\to\infty$,
\[
P(\max_{1\le j\le r_n}Y_{nj}^2\le \beta_n+\alpha_nx)=\prod^{r_n}_{j=1}(1-z_{nj})
\]
and
\[
P(\max_{1\le j\le p_n}Y_{nj}^2\le \beta_n+
\alpha_nx)=\prod^{p_n}_{j=1}(1-z_{nj}).
\]
By using Lemma~\ref{prod2sum} and \eqref{cdfG-rn} we have $\sum^{r_n}_{j=1}z_{nj}\to z=-\log G(x)$ which, together with Lemma~\ref{little},
implies $\sum^{p_n}_{j=1}z_{nj}\to z=-\log G(x)$. Once again we have from  Lemma~\ref{prod2sum}  that
 \[
\lim_{n\to\infty}P(\max_{1\le j\le p_n}Y_{nj}^2\le \beta_n+\alpha_nx)=\prod^{p_n}_{j=1}(1-z_{nj})=e^{-z}=G(x).
\]
This completes the proof of the lemma. \hfill$\blacksquare$

\begin{lemma}\label{bird} Let $Z_n$ be nonnegative random variables such that $(Z_n^2-\beta_n)/\alpha_n$ converges weakly to a cdf $G(x)$, where
$\alpha_n>0$ and $\beta_n>0$ are constants satisfying that $\lim_{n\to\infty}\alpha_n/\beta_n=0$. Then
\begin{equation}\label{ZZZ}
\frac{Z_n-\beta_n^{1/2}}{\alpha_n/(2\beta_n^{1/2})} \mbox{ converges weakly to } G.
\end{equation}
\end{lemma}

\noindent{\bf Proof.}  Set $W_n=(Z_n^2-\beta_n)/\alpha_n$. We have $Z_n^2=\beta_n+\alpha_nW_n=\beta_n(1+\frac{\alpha_n}{\beta_n}W_n)$. Then by Taylor's expansion
\[
Z_n=\beta_n^{1/2}(1+\frac{\alpha_n}{\beta_n}W_n)^{1/2}=\beta_n^{1/2}(1+\frac{\alpha_n}{2\beta_n}W_n+O_p(\frac{\alpha_n}{\beta_n})^2)
\]
and thus we have
\[
\frac{Z_n-\beta_n^{1/2}}{\alpha_n/(2\beta_n^{1/2})}=W_n+O_p(\frac{\alpha_n}{\beta_n}),
\]
which implies \eqref{ZZZ}. \hfill$\blacksquare$

\begin{lemma}\label{seed} (Lemma 2.2 of Jiang and Qi (2017))
Let $\{j_n,  n\ge 1\}$ and $\{x_n, n\geq 1\}$ be positive numbers with
$\lim_{n\to\infty}x_n=\infty$
and  $\lim_{n\to\infty}j_n x_n^{-1/2}(\log x_n)^{1/2}=\infty$.
For fixed $y\in \mathbb{R}$,  if
$\{c_{n,j}, 1\le j\le j_n, n\ge 1\}$ are real numbers
such that $\lim_{n\to\infty}\max_{1\le j\le j_n}|c_{n,j}x_n^{1/2}-1|=0$, then
\bea
& & \lim_{n\to\infty}\sum_{j=1}^{j_n}\big(1-\Phi((j-1)c_{n,j}+a(x_n)+b(x_n)y)\big)=e^{-y};\label{normal-app1}\\
& & \lim_{n\to\infty}\sum_{j=1}^{j_n}\frac{1}{(j-1)c_{n,j}+a(x_n)+b(x_n)y}\phi((j-1)c_{n,j}+a(x_n)+b(x_n)y\big)=e^{-y}.
\label{normal-app2}  \ \ \ \ \ \ \ \ \ \
\eea
\end{lemma}

\begin{lemma}\label{reiss} (Lemma 2.3 of Jiang and Qi (2017)  or Proposition 2.10 of Reiss (1981))
Let $\cal{B}$ be the collection of all Borel sets on $\mathbb{R}$.
  Then there exists a constant $C>0$ such that for all $r>k\geq 1$,
\beaa
& & \sup_{B\in \cal{B}}\Big|P\Big(\frac{r^{3/2}}{\sqrt{(r-k)k}}\big(U_{r-k+1:r}-\frac{r-k}{r}\big)\in B\Big)-\int_B(1+l_1(t)+l_2(t))\phi(t)dt \Big|\\
&\le &
C\cdot\Big(\frac{r}{(r-k)k}\Big)^{3/2}
\eeaa
where
%, $\phi(x)=\frac{1}{\sqrt{2\pi}}e^{-\frac{x^2}{2}}$ is the standard normal
%density function,
for $i=1,2$, $l_i(t)$ is a polynomial in $t$ of degree $\le 3i$, depending on $r$ and $k$,  and all of its coefficients  are of order $O(\big(\frac{r}{(r-k)k}\big)^{i/2})$.
\end{lemma}

\begin{lemma}\label{expansion1}
Define $V_{p_n-j+1: n-j}$ as in \eqref{vj}.  Assume that $k_n=n-p_n\to\infty$ and $k_n/n\to 0$ as $n\to\infty$.    Then for any $\delta_n>0$
such that $\delta_n\to\infty$ and $\delta_n=o(k_n^{1/6})$
\begin{equation}\label{eq1}
P(V_{p_n-j+1: n-j}>x)=(1+o(1))(1-\Phi(x))
\end{equation}
uniformly over $0\le x\le \delta_n$, $1\le j\le p_n-k_n$ as $n\to\infty$.
\end{lemma}

\noindent{\bf Proof.} Set $\beta_{nj}(x)=\frac{p_n-j}{n-j}+\frac{((p_n-j)k_n)^{1/2}}{(n-j)^{3/2}}x=\frac{p_n-j}{n-j}(1+\frac{k_n^{1/2}}{(n-j)^{1/2}(p_n-j)^{1/2}}x)$. Then
$1-\beta_{nj}(x)=\frac{n-p_n}{n-j}-\frac{((p_n-j)k_n)^{1/2}}{(n-j)^{3/2}}x=\frac{k_n}{n-j}(1-\frac{(p_n-j)^{1/2}}{(n-j)^{1/2}k_n^{1/2}}x)$,
and the density function of  $V_{p_n-j+1: n-j}$  is given by
\bea
h_j(x)&=&\frac{((p_n-j)k_n)^{1/2}}{(n-j)^{3/2}}f_{p_n-j+1:n-j}(\beta_{nj}(x))\nonumber\\
&=&\frac{((p_n-j)k_n)^{1/2}}{(n-j)^{3/2}}\frac{(n-j)!}{(p_n-j)!(k_n-1)!}\beta_{nj}(x)^{p_n-j}(1-\beta_{nj}(x))^{k_n-1}\nonumber\\
&=&\frac{(p_n-j)^{p_n-j+1/2}k_n^{k_n+1/2}}{(n-j)^{n-j+1/2}}\frac{(n-j)!}{(p_n-j)!k_n!}(1-\frac{(p_n-j)^{1/2}}{(n-j)^{1/2}k_n^{1/2}}x)^{-1}
\label{stirlingprod}\\
&&~\times(1+\frac{k_n^{1/2}}{(n-j)^{1/2}(p_n-j)^{1/2}}x)^{p_n-j}
(1-\frac{(p_n-j)^{1/2}}{(n-j)^{1/2}k_n^{1/2}}x)^{k_n}.\label{taylorprod}
\eea
To estimate $h_j(x)$, we need Stirling's formula:
\begin{equation}\label{stirling}
j!=j^{j+1/2}e^{-j+\varepsilon(j)}\sqrt{2\pi}, \mbox{ where } \frac{1}{12j+1}<\varepsilon(j)<\frac{1}{12j}
\end{equation}
 and Taylor's expansion: $1-t=\exp(\log(1-t))=\exp(-t-\frac{1}{2}t^2+O(t^3))$ as $t\to 0$.
By applying Stirling's formula to $k_n!$, $(p_n-j)!$ and $(n-j)!$, the product in \eqref{stirlingprod} is equal to $\frac{1+o(1)}{\sqrt{2\pi}}(1+o(1))$ for $|x|=o(k_n^{1/6})$
uniformly over $1\le j\le p_n-k_n$ as $n\to\infty$.  By applying Taylor's expansion to $(1+\frac{k_n^{1/2}}{(n-j)^{1/2}(p_n-j)^{1/2}}x)^{p_n-j}$ and
$(1-\frac{(p_n-j)^{1/2}}{(n-j)^{1/2}k_n^{1/2}}x)^{k_n}$, the product in \eqref{taylorprod} is equal to
\beaa
&&\exp\Big(\frac{(p_n-j)k_n^{1/2}x}{(n-j)^{1/2}(p_n-j)^{1/2}}-\frac12\frac{(p_n-j)k_nx^2}{(n-j)(p_n-j)}
+O(\frac{(p_n-j)k_n^{3/2}|x^3|}{(n-j)^{3/2}(p_n-j)^{3/2}})\Big)\\
&&\times\exp\Big(-\frac{k_n(p_n-j)^{1/2}x}{(n-j)^{1/2}k_n^{1/2}}-\frac12\frac{k_n(p_n-j)x^2}{(n-j)k_n}+
O(\frac{k_n(p_n-j)^{3/2}|x|^3}{(n-j)^{3/2}k_n^{3/2}})\Big)\\
&=&\exp\Big(\frac{(p_n-j)^{1/2}k_n^{1/2}x}{(n-j)^{1/2}}-\frac12\frac{k_nx^2}{(n-j)}
+O(\frac{k_n^{3/2}|x^3|}{(n-j)^{3/2}(p_n-j)^{1/2}})\Big)\\
&&\times\exp\Big(-\frac{k_n^{1/2}(p_n-j)^{1/2}x}{(n-j)^{1/2}}-\frac12\frac{(p_n-j)x^2}{(n-j)}+
O(\frac{(p_n-j)^{3/2}|x|^3}{(n-j)^{3/2}k_n^{1/2}})\Big)\\
&=&\exp\Big(-\frac{x^2}2+O((\frac1{k_n^{1/2}}+\frac1{(p_n-j)^{1/2}})|x|^3)\Big)\\
&=&\exp\Big(-\frac{x^2}2+o(1)\Big)
\eeaa
for $|x|=o(k_n^{1/6})$
uniformly over $1\le j\le p_n-k_n$ as $n\to\infty$. Therefore, we have
\begin{equation}\label{exact}
h_j(x)=\frac{1+o(1)}{\sqrt{2\pi}}\exp(-\frac{1}{2}x^2),~~~|x|=o(k_n^{1/6})
\end{equation}
uniformly over $1\le j\le p_n-k_n$ as $n\to\infty$.  Next we will give an estimate of the upper bound of $h_j(x)$ for large $x$. Note that
$\beta_{nj}(x)<1$ if and only if $x<\frac{k_n^{1/2} (p_n-j)^{1/2}}{(n-j)^{1/2}}=O(k_n^{1/2})$ uniformly over $1\le j\le p_n-k_n$ and thus  $\frac{k_n^{1/2}x}{(n-j)^{1/2}(p_n-j)^{1/2}}\le O(\frac{k_n^{1/2}}{n^{1/2}})\to 0$ uniformly over $0<x<\frac{k_n^{1/2} (p_n-j)^{1/2}}{(n-j)^{1/2}}$, $1\le j\le p_n-k_n$ as $n\to\infty$. Now by applying Taylor's expansion to  $(1+\frac{k_n^{1/2}}{(n-j)^{1/2}(p_n-j)^{1/2}}x)^{p_n-j}$ and inequality \[
1-t\le \exp(-t-\frac{1}{2}t^2), ~~~t\in (0,1)
\]
to $(1-\frac{(p_n-j)^{1/2}}{(n-j)^{1/2}k_n^{1/2}}x)^{k_n}$ we get
\begin{equation}\label{control}
h_j(x)\le \frac{1+o(1)}{\sqrt{2\pi}}\exp(-\frac{1}{2}x^2),~~~0<x<\frac{k_n^{1/2} (p_n-j)^{1/2}}{(n-j)^{1/2}}
\end{equation}
uniformly over $1\le j\le p_n-k_n$ as $n\to\infty$.

Assume that $0\le x\le \delta_n$. From \eqref{control} we have
\beaa
P(V_{p_n-j+1: n-j}>x)&=&\int^{\frac{k_n^{1/2} (p_n-j)^{1/2}}{(n-j)^{1/2}}}_xh_j(t)dt\\
&\le& (1+o(1))\int^{\frac{k_n^{1/2} (p_n-j)^{1/2}}{(n-j)^{1/2}}}_x\phi(t)dt\\
&\le& (1+o(1))(1-\Phi(x))
\eeaa
uniformly over $0\le x\le \delta_n$, $1\le j\le p_n-k_n$ as $n\to\infty$. Therefore, to complete the proof of the lemma,
it suffices to show that
\begin{equation}\label{inequality2}
P(V_{p_n-j+1: n-j}>x)\ge (1+o(1))(1-\Phi(x))
\end{equation}
uniformly over $0\le x\le \delta_n$, $1\le j\le p_n-k_n$ as $n\to\infty$.

From \eqref{tailapp}, we see that
\[
\Phi(y)-\Phi(x)=(1+o(1))(1-\Phi(x)),  ~~~0\le x=o(y)
\]
uniformly if $y\to\infty$.  Since $(\delta_nk_n^{1/6})^{1/2}=o(k_n^{1/6})$,
by using \eqref{exact} we have
\beaa
P(V_{p_n-j+1: n-j}>x)&\ge &\int^{(\delta_nk_n^{1/6})^{1/2}}_xh_j(t)dt\\
&=& (1+o(1))\int^{(\delta_nk_n^{1/6})^{1/2}}_x\phi(t)dt\\
&=& (1+o(1))(\Phi((\delta_nk_n^{1/6})^{1/2}) -\Phi(x))\\
&=&(1+o(1))(1-\Phi(x)),
\eeaa
proving \eqref{inequality2}.  This completes the proof of the lemma. \hfill $\blacksquare$

\subsection{Proofs of the Theorems}

\noindent\textbf{Proof of Theorem~\ref{thm2}}. We need to prove
\begin{equation}\label{thm2eq1}
\frac{1}{B_n}(\max_{1\leq j \leq p_n}|z_j|-A_n)\td \Lambda,
\end{equation}
where $A_n=c_n+\frac{1}{2}(1-c_n^2)^{1/2}(n-1)^{-1/2}a_n$,  $B_n=\frac{1}{2}(1-c_n^2)^{1/2}(n-1)^{-1/2}b_n$,
\[
c_n=\Big(\frac{p_n-1}{n-1}\Big)^{1/2},~~ b_n=b\Big(\frac{nc_n^2}{1-c_n^2}\Big), ~~a_n=a\Big(\frac{nc_n^2}{1-c_n^2}\Big)
\]
with $a(x)=(\log x)^{1/2}-(\log x)^{-1/2}\log(\sqrt{2\pi}\log x)$ and $b(x)=(\log x)^{-1/2}$ for $x>3$.

Fix $x\in \mathbb{R}$ and set $t_n=t_n(x)=c_n^2+c_n(1-c_n^2)^{1/2}(n-1)^{-1/2}(a_n+b_nx)$. For each $n\ge 2$, define $z_{nj}$ as in \eqref{znj}, that is, $z_{nj}=1-F_{p_n-j+1:n-j}(t_n(x))$ for $1\le j\le p_n$.

Since $Y_{nj}^2$ and $U_{p_n-j+1: n-j}$ are identically distributed, we have
\begin{equation}\label{aprod}
P\Big(\max_{1\leq j \leq r}Y_{nj}^2\le t_n(x)\Big)=\prod^r_{j=1}P(Y_{nj}^2\le t_n(x))=
\prod^r_{j=1}(1-z_{nj})
\end{equation}
for any $1\le r\le p_n$.

\noindent{\it Part 1}. First we show \eqref{thm2eq1} under condition \eqref{c1}. We will prove that
\begin{equation}\label{red_face}
\lim_{n\to\infty}P\Big(\max_{1\leq j \leq p_n}Y_{nj}^2\le t_n(x)\Big)=\exp(-e^{-x}).
\end{equation}

Let $j_n=[p_n^{5/8}]$, the integer part of $p_n^{5/8}$. For $1\le j\le j_n$, define
\[
u_{nj}=\frac{(n-j)^{3/2}}{((p_n-j)k_n)^{1/2}}\left(t_n(x)-\frac{p_n-j}{n-j}\right).
\]
Meanwhile, we rewrite
\[
t_n(x)=\frac{p_n-1}{n-1}+\frac{((p_n-1)k_n)^{1/2}}{(n-1)^{3/2}}(a_n+b_nx).
\]
Then we see that uniformly over $1\le j\le j_n$,
\beaa
u_{nj}&= & \Big(\frac{p_n-1}{n-1}-\frac{p_n-j}{n-j}\Big)\cdot \frac{(n-j)^{3/2}}{((p_n-j)k_n)^{1/2}}\\
& & ~~~~~~~~~~~~~~~~~~~~~~~~~~~~~~~~~~~~~~ + \Big(\frac{n-j}{n-1}\Big)^{3/2}\cdot \Big(\frac{p_n-1}{p_n-j}\Big)^{1/2}(a_n+b_nx)\\
& = & \Big(\frac{p_n-j}{p_n-1}\Big)^{-1/2}\cdot \Big(\frac{n-j}{n-1}\Big)^{1/2}\cdot
\Big(\frac{n-p_n}{p_n-1}\Big)^{1/2}\cdot\Big(\frac{n-1}{n}\Big)^{-1/2}\cdot \frac{j-1}{n^{1/2}}\\
& & ~~~~~~~~~~~~~~~~~~~~~~~~~~~~~~~~~~~~~~ + \Big(\frac{n-j}{n-1}\Big)^{3/2}\cdot \Big(\frac{p_n-j}{p_n-1}\Big)^{-1/2}(a_n+b_nx).
\eeaa
Now, $\frac{n-p_n}{p_n-1}=\frac{1-c_n^2}{c_n^2}\sim \frac{n}{p_n}$. Also, given $\tau\in \mathbb{R}$, trivially $\Big(\frac{p_n-j}{p_n-1}\Big)^{\tau}=1+O(\frac{j-1}{p_n})$ and $\Big(\frac{n-j}{n-1}\Big)^{\tau}=1+O(\frac{j-1}{n})$ uniformly for all $1\leq j \leq j_n.$ Since $a_n \sim (\log p_n)^{1/2}$ and $b_n\sim (\log p_n)^{-1/2}$, we have
\begin{eqnarray}\label{onlyone}
u_{nj} &=& \frac{(1-c_n^2)^{1/2}}{n^{1/2}c_n}(j-1)(1+o(1)) + a_n+ b_nx + O\Big(\frac{(j-1)(\log p_n)^{1/2}}{p_n}\Big)\nonumber\\
& = & \frac{(1-c_n^2)^{1/2}}{n^{1/2}c_n}(j-1)(1+o(1)) + a_n+ b_nx
\end{eqnarray}
uniformly for all $1\leq j \leq j_n$ as $n\to\infty$.

In Lemma~\ref{reiss}, take $r=n-j$ and $k=n-p_n$ to have
\[
\sup_{B\in \cal{B}}\Big|P(V_{p_n-j+1:n-j}\in B)-\int_B(1+l_1(t)+l_2(t))\phi(t)dt\Big|=O(\frac{1}{p_n^{3/2}})
\]
uniformly over $1\le j\le j_n$ as $n\to\infty$,  where
\begin{equation}\label{vj}
V_{p_n-j+1:n-j}=\frac{(n-j)^{3/2}}{((p_n-j)(n-p_n))^{1/2}}\Big(U_{p_n-j+1:n-j}-\frac{p_n-j}{n-j}\Big)
\end{equation}
and
where,
for $i=1,2$, $l_i(t)$ is a polynomial in $t$ of degree $\le 3i$, depending on $n$,  and all of its coefficients  are of order $O((1/p_n)^{i/2})$ uniformly over $1\le j\le j_n$ as $n\to\infty$.  Now, by taking $B=(u_{nj}, \infty)$ we obtain
\begin{equation}\label{znjthm2}
z_{nj}=P(V_{p_n-j+1:n-j}>u_{nj})
=\int^\infty_{u_{nj}}(1+l_1(t)+l_2(t))\phi(t)dt+O(p_n^{-3/2})
\end{equation}
uniformly for $1\le j\le j_n$ as $n\to\infty$. From L'Hospital's rule, we have  that for any $r\ge 0$
\begin{equation}\label{tailapp}
\int^\infty_xt^r\phi(t)dt\sim x^{r-1}\phi(x)~~~~~\mbox{ as } x\to\infty.
\end{equation}
Since $\min_{1\le j\le j_n}u_{nj}\to\infty$ as $n\to\infty$ by (\ref{onlyone}), it  follows from \eqref{tailapp} that
\[
\int^\infty_{u_{nj}}t^r\phi(t)dt\sim (u_{nj})^{r-1}\phi(u_{nj})
\]
holds uniformly over $1\le j\le j_n$. Furthermore, since the coefficients of
$l_i(t)$ are uniformly bounded by $O((1/p_n)^{i/2})$ for $i=1,2$,
we have
\begin{eqnarray}\label{3terms}
&&\sum^{j_n}_{j=1}\int^\infty_{u_{nj}}(1+l_1(t)+l_2(t))\phi(t)dt\\
&=&(1+o(1))\sum^{j_n}_{j=1}\frac{\phi(u_{nj})}{u_{nj}}
+O(\frac{1}{p_n^{1/2}})\sum^{j_n}_{j=1}u_{nj}^3 \frac{\phi(u_{nj})}{u_{nj}}+
O(\frac{1}{p_n})\sum^{j_n}_{j=1}u_{nj}^6 \frac{\phi(u_{nj})}{u_{nj}}.\nonumber
\end{eqnarray}
In Lemma~\ref{seed}, by taking $x_n=nc_n^2/(1-c_n^2)$ and and define $c_{n,j}$ such that
$u_{nj}=(j-1)c_{n,j}+a(x_n)+b(x_n)x$ for $1\le j\le j_n$ with $c_{n,1}=x_n^{-1/2}$. It follows from \eqref{onlyone} that
$c_{n,j}=x_n^{-1/2}(1+o(1))$ uniformly over $1\le j\le j_n$ as $n\to\infty$, which implies $\lim_{n\to\infty}\max_{1\le j\le j_n}|c_{n,j}x_n^{1/2}-1|=0$. Then we get
\begin{equation}\label{sum1}
\lim_{n\to\infty}\sum^{j_n}_{j=1}\frac{\phi(u_{nj})}{u_{nj}}= e^{-x}.
\end{equation}
 We will show that the second term and the third term on the line below \eqref{3terms} converge to zero as $n\to\infty$. By noting that $\frac{(1-c_n^2)^{1/2}}{n^{1/2}c_n}\sim p_n^{-1/2}$ we have
\begin{equation}\label{unjbound}
u_{nj}^3=O(\frac{j_n^3}{p_n^{3/2}})=O(p_n^{3/8})
\end{equation}
uniformly over $1\le j\le j_n$. Thus, it follows from \eqref{normal-app2} that
\[
O(\frac{1}{p_n^{1/2}})\sum^{j_n}_{j=1}u_{nj}^3 \frac{\phi(u_{nj})}{u_{nj}}=O(\frac{1}{p_n^{1/8}})\sum^{j_n}_{j=1}\frac{\phi(u_{nj})}{u_{nj}}=O(\frac{1}{p_n^{1/8}})\to 0
\]
as $n\to\infty$. Similarly, we have
\[
O(\frac{1}{p_n})\sum^{j_n}_{j=1}u_{nj}^6 \frac{\phi(u_{nj})}{u_{nj}}=O(\frac{1}{p_n^{1/4}})\to 0
\]
as $n\to\infty$. Therefore, by combining \eqref{3terms}, \eqref{znjthm2} and \eqref{sum1} we get
\begin{equation}\label{sum2}
\sum^{j_n}_{j=1}z_{nj}\to e^{-x} ~~ \mbox{ as }n\to\infty.
\end{equation}
It follows from \eqref{onlyone} that
\[
u_{nj_n}\ge \frac{(j_n-1)(1+o(1))}{p_n^{1/2}}\ge \frac{p_n^{1/8}}{2}
\]
for all large $n$. Then we estimate $z_{nj_n}$ by using \eqref{znjthm2} with $j=j_n$ and \eqref{unjbound}
\beaa
z_{nj_n}&=&\big(\frac{1+o(1)}{u_{nj_n}}+O(\frac{u_{nj_n}^2}{p_n^{1/2}})+O(\frac{u_{nj_n}^5}{p_n})\big)\phi(u_{nj_n})+
O(\frac{1}{p_n^{3/2}})\\
&=&O(\exp(-\frac12u_{nj_n}^2))+O(\frac{1}{p_n^{3/2}})\\
&=&O(\exp(-\frac18p_n^{1/4}))+O(\frac{1}{p_n^{3/2}})\\
&=&O(\frac{1}{p_n^{3/2}})
\eeaa
as $n\to\infty$. From \eqref{monotonez} we have $\sum^{p_n}_{j=j_n+1}z_{nj}\le p_nz_{nj_n}=O(p_n^{-1/2})\to 0$
as $n\to\infty$, which together with \eqref{sum2} yields
\[
\sum^{p_n}_{j=1}z_{nj}\to e^{-x}~~~\mbox{ as }n\to\infty.
\]
We can also prove from \eqref{znjthm2} that $z_{n1}\to 0$ as $n\to\infty$.  In view of  \eqref{aprod} and Lemma~\ref{prod2sum} we conclude
\eqref{red_face}, ie.,
\[
\frac{\max_{1\le j\le p_n}Y_{nj}^2-\beta_n}{\alpha_n}\td \Lambda,
\]
where $\alpha_n=c_n(1-c_n^2)^{1/2}(n-1)^{-1/2}b_n$  and $\beta_n=c_n^2+c_n(1-c_n^2)^{1/2}(n-1)^{-1/2}a_n$.
Since
\[
\frac{\alpha_n}{\beta_n}\le \frac{c_n(1-c_n^2)^{1/2}(n-1)^{-1/2}b_n}{c_n(1-c_n^2)^{1/2}(n-1)^{-1/2}a_n}=\frac{b_n}{a_n}\to 0
\]
as $n\to\infty$, we can apply Lemma~\ref{bird} and get that
\begin{equation}\label{love}
\Lambda_n:=\frac{\max_{1\le j\le p_n}Y_{nj}-\beta_n^{1/2}}{\alpha_n/(2\beta_n^{1/2})}\td \Lambda.
\end{equation}
Recall that  $A_n=c_n+\frac{1}{2}(1-c_n^2)^{1/2}(n-1)^{-1/2}a_n=c_n(1+o(1))$ and
$B_n=\frac{1}{2}(1-c_n^2)^{1/2}(n-1)^{-1/2}b_n\sim \frac12(n-1)^{-1/2}(\log p_n)^{-1/2}$.
Then
\beaa
\beta_n^{1/2}&=&c_n(1+\frac{(1-c_n^2)^{1/2}a_n}{(n-1)^{1/2}c_n})^{1/2}\\
&=&c_n(1+\frac{(1-c_n^2)^{1/2}a_n}{2(n-1)^{1/2}c_n}+O(\frac{(1-c_n^2)a_n^2}{(n-1)c_n^2}))\\
&=&A_n+O(\frac{\log p_n}{(np_n)^{1/2}})\\
&=&A_n+o(B_n)
\eeaa
and
\beaa
\frac{\alpha_n}{2\beta_n^{1/2}}&=&B_n(1+o(1)),
\eeaa
which, together with \eqref{love}, yield
\beaa
\frac{\max_{1\le j\le p_n}Y_{nj}-A_n}{B_n}&=&\frac{\Lambda_n\alpha_n/(2\beta_n^{1/2})+\beta_n^{1/2}-A_n}{B_n}\\
&=&\frac{\Lambda_nB_n(1+o(1))+o(B_n)}{B_n}\\
&=&(1+o(1))\Lambda_n+o(1)\\
&\td& \Lambda
\eeaa
ie., \eqref{thm2eq1} holds. The proof of {\it Part 1} is completed.

\noindent{\it Part 2.} We will show \eqref{thm2eq1} under condition \eqref{c2}.  First, it follows from condition \eqref{c2} that
\[
\frac{nc_n^2}{1-c_n^2}=\frac{n(p_n-1)}{n-p_n}\sim \frac{n^2}{k_n}\to\infty
\]
as $n\to\infty$. Noting that $n\le n^2/k_n\le n^2$, we get that
\[
a_n=a(\frac{nc_n^2}{1-c_n^2})\sim \big(\log(\frac{nc_n^2}{1-c_n^2})\big)^{1/2}
\]
is of order $(\log n)^{1/2}$ as $n\to\infty$.

Use the same notation as in Part 1. Recall that $p_n/n\to 1$, $k_n=n-p_n=o(n)$ and $\log n=o(k_n^{1/3})$ as $n\to\infty$.
In order to use both Lemmas~\ref{seed} and \ref{expansion1}, we take $x_n=\frac{nc_n^2}{1-c_n^2}$. Define $j_n=[5(\log n)^{1/2}\sqrt{x_n}]+1$. Then $j_n\sim \frac{5n(\log n)^{1/2}}{\sqrt{k_n}}=o(n)$, which implies $1\le j_n\le p_n-k_n$ for all large $n$.
Define $c_{n,j}$ for $1\le j\le j_n$ in the same way as in {\it Part 1}.
Similar to the proof of
\eqref{onlyone} we can show that
\begin{equation}\label{onlytwo}
u_{nj} =\frac{(1-c_n^2)^{1/2}}{n^{1/2}c_n}(j-1)(1+o(1)) + a_n+ b_nx
\end{equation}
uniformly for $1\leq j \leq j_n$ as $n\to\infty$. Then $u_{nj}=O((\log n)^{1/2})=o(k_n^{1/6})$ uniformly for $1\leq j \leq j_n$ as $n\to\infty$. We can also verify that all conditions in Lemma~\ref{seed} are satisfied.
Thus, from
\eqref{eq1} and \eqref{normal-app1} we have
\begin{equation}\label{limit1}
\sum^{j_n}_{j=1}z_{nj}=(1+o(1))\sum^{j_n}_{j=1}(1-\Phi(u_{nj}))\to e^{-x}
\end{equation}
as $n\to\infty$.

Next, we will show that
\begin{equation}\label{limit2}
\lim_{n\to\infty}\sum^{p_n}_{j=j_n+1}z_{nj}=0.
\end{equation}
Note that \eqref{eq1} holds uniformly over $1\le j\le p_n-k_n$ and
$u_{nj_n}\ge 4(\log n)^{1/2}$ for all large $n$. By employing \eqref{eq1} with $j=j_n$ and $x=4(\log n)^{1/2}$ and
using equation \eqref{monotonez} and Lemma~\ref{expansion1} we have
\beaa
\sum^{p_n}_{j=j_n+1}z_{nj}&\le &nz_{nj_n}\\
&\le & nP(V_{p_n-j_n+1:n-j_n}>u_{nj_n})\\
&\le & nP(V_{p_n-j_n+1:n-j_n}\ge 4(\log n)^{1/2})\\
&=&(1+o(1))n(1-\Phi( 4(\log n)^{1/2}))\\
&=&O(\frac1n)\\
&\to& 0.
\eeaa
Thus, we obtain that $\sum^{p_n}_{j=1}z_{nj}\to e^{-x}$ for any $x$.  Then equation \eqref{red_face} follows from  equation \eqref{aprod} and Lemma~\ref{prod2sum}.  The rest of the proof will follow from the same lines in the proof of the first part. Again Lemma~\ref{bird} will be used.  The details are omitted.
\hfill $\blacksquare$\\

%so that $4(\log n)^{1/2}\le u_{nj}=O((\log n)^{1/2})=o(k_n^{1/6})$ uniformly over $1\le j\le j_n$.

\noindent{\bf Proof of Theorem~\ref{thm3}.}  Recall that $a_n$ is given by
\begin{equation}\label{an}
\frac{1}{(k_n-1)!}\int^{a_n}_0t^{k_n-1}e^{-t}dt=\frac{k_n}{n}.
\end{equation}

By using integration by parts, we have
\[
y^{k_n}e^{-y}=k_n\int^y_0t^{k_n-1}e^{-t}dt-\int^y_0t^{k_n}e^{-t}dt\ge k_n\int^y_0t^{k_n-1}e^{-t}dt-y\int^y_0t^{k_n-1}e^{-t}dt,
\]
which implies
 \begin{equation}\label{integral}
\frac{1}{k_n}y^{k_n}e^{-y}\le \int^y_0t^{k_n-1}e^{-t}dt\le \frac{1}{k_n-y}y^{k_n}e^{-y}, ~~~0\le y<k_n.
 \end{equation}
By using Stirling's formula \eqref{stirling}, we have under condition~\eqref{c3} that
\[
y_n:=(\frac{k_n!k_n}{nk_n^{k_n}})^{1/k_n}
=\exp(\frac32\frac{\log k_n}{k_n}-1+\frac{\varepsilon(k_n)+\log\sqrt{2\pi}}{k_n})\exp(-\frac{\log n}{k_n})
\to 0
\]
Since $te^{-t}$ is strictly increasing in $(0,1)$, for all large $n$ such that $y_n<1$  define
$\varepsilon_n$ as the unique solution to $te^{-t}=y_n$ in $(0,1)$, that is, $\varepsilon_ne^{-\varepsilon_n}=y_n$, which implies that
$\varepsilon_n\to 0$ and $\varepsilon_n\sim y_n$ as $n\to\infty$ and
\[
\frac{1}{k_n!}(k_n\varepsilon_n)^{k_n}e^{-k_n\varepsilon_n}=\frac{k_n}{n}
\]
for all large $n$. Then it follows from the the first inequality in \eqref{integral} that
\[
\frac{1}{(k_n-1)!}\int^{k_n\varepsilon_n}_0t^{k_n-1}e^{-t}dt\ge \frac{k_n}{n}
\]
for all large $n$, which together with \eqref{an} implies that $a_n\le k_n\varepsilon_n$ for all large $n$,and thus $a_n=o(k_n)$ as $n\to\infty$. By plugging $y=a_n$ in \eqref{integral} and using  $a_n\le k_n\varepsilon_n$ for large $n$  we conclude
\[
\int^{a_n}_0t^{k_n-1}e^{-t}dt=\frac{1}{k_n}a_n^{k_n}e^{-a_n}(1+O(\varepsilon_n))
\]
which implies
\[
\frac{k_n}{n}=\frac{1}{(k_n-1)!}\int^{a_n}_0t^{k_n-1}e^{-t}dt=\frac{1}{k_n!}a_n^{k_n}e^{-a_n}(1+O(\varepsilon_n))
\]
as $n\to\infty$, and consequently
\begin{equation}\label{constant}
\frac{n}{k_n!k_n}a_n^{k_n}e^{-a_n}=1+o(1)~~~\mbox{as }n\to\infty.
\end{equation}

Define $z_{nj}$ as in \eqref{znj} with $t_n=t_n(x)=1-\frac{a_n}{n}(1-\frac{x}{k_n})$ for any fixed $x$. Then  $1-t_n=\frac{a_n}{n}(1-\frac{x}{k_n})=o(\frac{k_n}{n})=o(\frac{\log n}{n})$ as $n\to\infty$, where we have used the fact that
$k_n=n-p_n\to\infty$ and $k_n=o(\log n)$ from \eqref{c3}.
This implies $n(1-t_n)^2\to 0$ as $n\to\infty$.

It is easy to verify the following expression
\[
1-t=e^{-t-d(t)}, ~~~0\le t\le 1/2
\]
where $0\le d(t)\le t^2$ for $0\le t\le 1/2$. Then
\beaa
\int^1_{t_n}t^{p_n-j}(1-t)^{k_n-1}dt
&=&\int^{1-t_n}_0(1-t)^{p_n-j}t^{k_n-1}dt\\
&=&(1+o(1))\int^{1-t_n}_0t^{k_n-1}e^{-(p_n-j)t}dt\\
&=&(1+o(1))\int^{1-t_n}_0t^{k_n-1}e^{-(n-j)t}e^{k_nt}dt\\
&=&(1+o(1))\int^{1-t_n}_0t^{k_n-1}e^{-(n-j)t}dt\\
&=&\frac{1+o(1)}{(n-j)^{k_n}}\int^{(1-t_n)(n-j)}_0t^{k_n-1}e^{-t}dt.
\eeaa
Furthermore, by using Stirling's formula \eqref{stirling}, we get
\beaa
\frac{(n-j)!}{(p_n-j)!}
&=&\frac{(n-j)^{n-j+1/2}e^{-(n-j)+\varepsilon(n-j)}}{(p_n-j)^{p_n-j+1/2}e^{-(p_n-j)+\varepsilon(p_n-j)}}\\
&=&(1+\frac{k_n}{p_n-j})^{p_n-j+1/2}(n-j)^{k_n}e^{-k_n}e^{\varepsilon(p_n-j)-\varepsilon(n-j)}\\
&=&(1+o(1))(n-j)^{k_n}
\eeaa
uniformly over $1\le j\le j_n$, where $j_n:=p_n-k_n^3$. Then from \eqref{znj=} we obtain that
\begin{equation}\label{znj-est}
z_{nj}=\frac{1+o(1)}{(k_n-1)!}\int^{(1-t_n)(n-j)}_0t^{k_n-1}e^{-t}dt
\end{equation}
uniformly over $1\le j\le j_n$.

Note that $n(1-t_n)=o(k_n)$. Then it follows from \eqref{znj-est} and \eqref{integral} that
\beaa
z_{nj}&=&\frac{1+o(1)}{k_n!}(n-j)^{k_n}(1-t_n)^{k_n}e^{-(n-j)(1-t_n)}\\
&=&\frac{1+o(1)}{1-t_n}\int^{(n-j+1)(1-t_n)}_{(n-j)(1-t_n)}t^{k_n}e^{-t}dt
\eeaa
uniformly over $1\le j\le j_n$, and thus
\beaa
\sum^{j_n}_{j=1}z_{nj}&=&\frac{1+o(1)}{k_n!}\sum^{j_n}_{j=1}(n-j)^{k_n}(1-t_n)^{k_n}e^{-(n-j)(1-t_n)}\\
&=&\frac{1+o(1)}{k_n!}\sum^{j_n}_{j=1}\frac{1}{1-t_n}\int^{(n-j+1)(1-t_n)}_{(n-j)(1-t_n)}t^{k_n}e^{-t}dt\\
&=&\frac{1+o(1)}{k_n!}\frac{1}{1-t_n}\int^{n(1-t_n)}_{(n-j_n)(1-t_n)}t^{k_n}e^{-t}dt\\
&=&\frac{1+o(1)}{k_n!}\frac{1}{1-t_n}\big(\int^{n(1-t_n)}_0t^{k_n}e^{-t}dt-\int^{k_n^3(1-t_n)}_0t^{k_n}e^{-t}dt\big).
\eeaa
The second integral above is dominated by the first one  since $t^{k_n}e^{-t}$ is increasing over $(0,k_n)$ and
$k_n^3(1-t_n)/(n(1-t_n))=k_n^3/n\to 0$ as $n\to\infty$. Therefore,  in view of \eqref{integral} and \eqref{constant} we have
\beaa
\sum^{j_n}_{j=1}z_{nj}&=&\frac{1+o(1)}{k_n!}\frac{1}{1-t_n}\int^{n(1-t_n)}_0t^{k_n}e^{-t}dt\nonumber\\
&=&\frac{1+o(1)}{(k_n+1)!}\frac{1}{1-t_n}(n(1-t_n))^{k_n+1}e^{-n(1-t_n)}\nonumber\\
&=&\frac{(1+o(1))n}{(k_n+1)!}(n(1-t_n))^{k_n}e^{-n(1-t_n)}\\
&=&\frac{(1+o(1))n}{k_n!k_n}(n(1-t_n))^{k_n}e^{-n(1-t_n)}\\
&=&\frac{(1+o(1))n}{k_n!k_n}a_n^{k_n}(1-\frac{x}{k_n})^{k_n}e^{-a_n+O(\frac{a_n}{k_n})}\\
&=&\frac{(1+o(1))n}{k_n!k_n}a_n^{k_n}e^{-a_n}e^{-x}\\
&\to& e^{-x}
\eeaa
as $n\to\infty$. Therefore, it follows from Lemma~\ref{little} that  $\sum^{p_n}_{j=1}z_{nj}\to e^{-x}$ as $n\to\infty$. It is easy to conclude that $z_{n1}\to 0$ as $n\to\infty$ from \eqref{znj-est}, \eqref{integral} and the above estimates. Accordingly, by taking $\beta_n=1-\frac{a_n}{n}$ and $\alpha_n=\frac{a_n}{nk_n}$ with $G(x)=\Lambda(x)$ in Lemmas~\ref{less} and \ref{bird} we conclude
that
\beaa
\frac{\max_{1\le j\le p_n}|z_j|-(1-a_n)^{1/2}}{a_n/(2nk_n)}&=&\frac{1}{(1-a_n)^{1/2}}\frac{\max_{1\le j\le p_n}|z_j|-(1-a_n)^{1/2}}{a_n/(nk_n2(1-a_n)^{1/2})}\\
&=&(1+o(1))\frac{\max_{1\le j\le p_n}|z_j|-(1-a_n)^{1/2}}{a_n/(nk_n2(1-a_n)^{1/2})}\\
&\td& \Lambda
\eeaa
as $n\to\infty$.
This completes the proof of the theorem. \hfill $\blacksquare$

\noindent{\bf Proof of Theorem~\ref{thm4}.}
We first show that
\begin{equation}\label{ss1}
\frac{n^{1+1/k}}{((k+1)!)^{1/k}}(\max_{1\le j\le p_n}|z_j|^2-1)\td W_{k}
\end{equation}

Fix $x<0$. Let $t_n=t_n(x)=1+\frac{((k+1)!)^{1/k}}{n^{1+1/k}}x$. Then $t_n\in (0,1)$ for all large $n$.
Since $n(1-t_n)\to 0$ as $n\to\infty$, we have
\[
(1-t)^{p_n-j}=1+o(1)~~~~\mbox{uniformly over $0\le  t\le 1-t_n$,  $1\le j\le p_n$}.
\]
Therefore, we have from \eqref{znj} and \eqref{znj=} that
\beaa
z_{nj}=1-F_{p_n-j+1:n-j}(t_n)&=&\frac{(n-j)!}{(p_n-j)!(k-1)!}\int^{1-t_n}_0(1-t)^{p_n-j}t^{k-1}dt\\
&=&\frac{(n-j)!}{(p_n-j)!(k-1)!}(1+o(1))\int^{1-t_n}_0t^{k-1}dt\\
&=&\frac{(n-j)!}{(p_n-j)!k!}(1+o(1))(1-t_n)^k\\
&=&\frac{(n-j)!}{(p_n-j)!}(1+o(1))\frac{(k+1)(-x)^{k}}{n^{k+1}}
\eeaa
uniformly over $1\le j\le p_n=n-k$. Since
\[
\frac{(n-j)!}{(p_n-j)!}\le n^k
\]
we have $\max_{1\le j\le p_n}z_{nj}=O(1/n)\to 0$ as $n\to\infty$. To complete the proof of \eqref{ss1},
by using \eqref{prod-sum} we need to show that
\begin{equation}\label{ss2}
\lim_{n\to\infty}\sum^{p_n}_{j=1}\frac{(n-j)!}{(p_n-j)!}\frac{k+1}{n^{k+1}}=1.
\end{equation}
Let $\{j_n\}$ be a sequence of integers such that $j_n\to\infty$ and $j_n/n\to 0$ as $n\to\infty$. Then
\[
\sum^{p_n}_{j=n-j_n+1}\frac{(n-j)!}{(p_n-j)!}\frac{k+1}{n^{k+1}}\le \frac{(k+1)j_n}{n}\to 0
\]
as $n\to\infty$, and
\[
\frac{(n-j)!}{(p_n-j)!}=(p_n-j+1)^k\prod^k_{\ell=1}(1+\frac{\ell-1}{p_n-j+1})=(p_n-j+1)^k(1+o(1))
\]
uniformly over $1\le j\le n-j_n$, which implies that
\beaa
\sum^{p_n-j_n}_{j=1}\frac{(n-j)!}{(p_n-j)!}\frac{k+1}{n^{k+1}}&=&\frac{(1+o(1))(k+1)}{n}\sum^{p_n-j_n}_{j=1}(\frac{p_n-j+1}{n})^k\\
&=&\frac{(1+o(1))(k+1)}{n}\sum^{p_n}_{j=j_n+1}(\frac{j}{n})^k\\
&\to& (k+1)
\int^1_{0}t^kdt=1
\eeaa
as $n\to\infty$. This proves \eqref{ss2} and thus we obtain \eqref{ss1}.

Finally, the theorem follows from Lemma~\ref{bird} with
$\alpha_n=\frac{((k+1)!)^{1/k}}{n^{1+1/k}}$ and $\beta_n=1$.
This completes the proof.
\hfill$\blacksquare$

\vspace{20pt}

\noindent\textbf{Acknowledgements} We would like to thank an anonymous referee for his/her careful reading of the original version of the paper and pointing out some imperfections in the proofs. Gui's work was partially supported by the program for the Fundamental Research Funds for the Central Universities (2014RC042).

\baselineskip 12pt
\def\ref{\par\noindent\hangindent 25pt}


\begin{thebibliography}{AA} %{SOSL90}

%\bibitem{Abramowitz1972} Abramowitz, M. and Stegun, I. A. (1972). Handbook of Mathematical Functions. {\it Dover, New York}.

%\bibitem{Nobel}
%Akemann, G., Baik, J. and Francesco, P. D. (2001).  The Oxford Handbook of Random Matrix Theory (Oxford Handbooks in Mathematics). {\it Oxford %University Press}.

%\bibitem{Akemann}
%Akemann, G. and Bender, M.  (2010). Interpolation between Airy and Poisson statistics for unitary chiral
%non-Hermitian random matrix ensembles. {\it J. Math. Phys.} 51(10): 103524.
%{\red The above one is uncited}

%\bibitem{Akemann2012}
%Akemann, G. and Burda, Z. (2012). Universal microscopic correlation functions for products of
%independent Ginibre matrices. {\it J. Phys. A: Math. Theor.} 45(46), 465201.

%\bibitem{Bai}
%Bai, Z. D. (1999). Methodologies in spectral analysis of large dimensional random matrices,
%a review. {\it Statistica Sinica} 9, 9611-9677.

%\bibitem{BaiYin}
%Bai, Z. D., Yin, Y. Q. and Krishnaiah, P. R. (1987). On the limiting empirical distribution
%function of the eigenvalues of a multivariate F matrix. {\it Theory Probab. Appl.} 32, 490-500.

\bibitem{Baik}
Baik, J, Deift, P. and Johansson, K. (1999). On the distribution of the length of the longest increasing
subsequence of random permutations. {\it Jour. Amer. Math. Soc.} 12(4), 1119-1178.

\bibitem{bala}
Balakrishnan, N. and Cohen, A. C. (1991).
Order Statistics and Inference: Estimation Methods. {\it Academic Press}.


\bibitem{BS}
Bl\"umel, R. and Smilansky, U. (1988). Classical irregular scattering and its quantum-mechanical implications.
Phys. Rev. Lett. 60, 477-480.



%\bibitem{Bor}
%Bordenave, C. (2011). On the spectrum of sum and product of non-Hemitian random matrices. {\it Elect. Comm. in Probab.} 16, 104-113.

%\bibitem{Bordoin}
%Borodin, A., Okounkov, A. and Olshanski, G. (2000).  Asymptotics of Plancherel
%measures for symmetric groups. {\it J. Amer. Math. Soc.} 13(3), 481-515.
%{\red The above one is uncited}

%\bibitem{braaksma}
%Braaksma, B. L. J. (1962-1964).
%Asymptotic expansions and analytic continuations for a class of Barnes-integrals. {\it Compositio Mathematica} 15, 239-341.
%{\red The above one is uncited}

%\bibitem{B}
%Burda, Z. (2013).  Free products of large random matrices - a short review of recent developments. {\it J. Phys. Conf. Ser. 473, 012002}. Also %available at {\tt http://arxiv.org/pdf/1309.2568v2.pdf}.

%\bibitem{BJW}
%Burda, Z., Janik, R. A. and Waclaw, B. (2010).  Spectrum of the product of independent random Gaussian matrices.  {\it Phys. Rev.} E 81, %041132.


%\bibitem{DJ69}
%Dharmadhikari, S.W. and Jogdeo, K. (1969). Bounds on
%moments of certain random variables. {\em Ann. Math. Statist.}  {\bf
%40},  1506-1508.

%\bibitem{Cha}
%Chafa\"{\i}, D. and P\'{e}ch\'{e}, S. (2014).
%A note on the second order universality at the edge of Coulomb gases on the plane. {\it Journal of Statistical Physics}
%156(2), 368-383.

%\bibitem{Chow}
%Chow, Y. S. and Teicher, H. (2003). Probability Theory: Independence, Interchangeability, Martingales. {\it Springer}, 3rd edition.

%\bibitem{Collins}
%Collins, B. (2005). Product of random projections, Jacobi ensembles and universality
%problems arising from free probability. {\it Probab. Theory Relat. Fields} 133, 315-344.

\bibitem{CD}
Couillet, R.  and Debbah, M. (2011). Random matrix methods for wireless communications. Cambridge
Univ Press, 2011.


%\bibitem{Dembo}
%Dembo, A. and Zeitouni, O. (1998). Large Deviations Techniques and Applications.
%{\it Springer}, 2nd edition.

%\bibitem{DI}
%Di Francesco, P.,  Gaudin, F. M., Itzykson, C. and Lesage, F. (1994). Laughlin's wave functions, Coulomb gases and expansions of the %discriminant.  {\it International Journal of Modern Physics A} 9, 4257-4351.


\bibitem{Evans}
Diaconis, P. and Evans, S. (2001). Linear functionals of eigenvalues of random matrices. {\it Transactions Amer. Math. Soc.} 353, 2615-2633.

\bibitem{Dong}
Dong, Z., Jiang, T. and Li, D. (2012).  Circular law and arc law for truncation of random unitary matrix.  {\it Journal of Mathematical Physics}  53, 013301-14.

%\bibitem{Dyson}
%Dyson F.J. (1962).  Statistical theory of the energy levels of complex systems. I {\it J. Math. Phys.} 3, 140-156.


%\bibitem{Eaton}
%Eaton, M. (2007). Multivariate Statistics: A Vector Space Approach. {\it IMS Lecture Notes Monograph}  53.

%\bibitem{Edelman}
%Edelman, A. (1997). The probability that a random real Gaussian matrix has k real eigenvalues,
%related distributions, and the circular law. {\it J. Multivariate Anal.} 60, 203-232.

\bibitem{Erdos}
Erd\H{o}s, L., Knowles, A., Yau, H. and Yin, J. (2012). Spectral statistics of Erd\H{o}s-R\'{e}nyi graphs II: eigenvalue spacing and the extreme eigenvalues. {\it Comm. Math. Phy.} 314(3), 587-640.


\bibitem{For}
Forrester, P.J. (2010). {\it Log-gases and random matrices}. Number 34. Princeton Univ Press.

%\bibitem{ForresterM}
%Forrester, P. J. and Mays, A. (2011). Pfaffian point process for the Gaussian real generalised eigenvalue problem. {\it Probability Theory and %Related Fields} 154(1-2), 1-47.

%\bibitem{ForresterN}
%Forrester, P. J. and Nagao, T. (2008). Skew orthogonal polynomials and the partly symmetric real Ginibre ensemble.  {\it J. Phys. A: Math. %Theor.} 41, 375003.

\bibitem{Ginibre}
Ginibre, J. (1965). Statistical ensembles of complex, quaternion, and real matrices. {\it J. Math. Phys.}
6, 440-449.

%\bibitem{Goetz}
%G\"{o}tze, F. and Tikhomirov, T. (2010). On the asymptotic spectrum of products of independent random
%matrices. {\tt http://arxiv.org/pdf/1012.2710v3.pdf}.

%\bibitem{Haake}
%Haake, F. (2010). Dissipative systems. In ``Quantum Signatures of Chaos" {\it Springer Series in Synergetics} 54, 279-339.

\bibitem{HF}
Hiai, F. and Petz, D. (2000). The Semicircle Law, Free Random Variables and Entropy,
Mathematical Surveys and Monographs, Vol. 77, American Mathematical Society.


%\bibitem{Hough}
%Hough, J. B., Krishnapur, M., Peres, Y. and Vir\'{a}g, B. (2009). Zeros of Gaussian Analytic Functions and Determinantal Point Processes. {\it %American Mathematical Society.}

\bibitem{Jiang09}
Jiang, T.  (2009). Approximation of Haar distributed matrices and limiting distributions of eigenvalues of Jacobi ensembles.  {\it Probability Theory and Related Fields} 144(1), 221-246.

\bibitem{Jiang10}
Jiang, T. (2010). The entries of Haar-invariant matrices from the classical compact groups.  {\it Journal of Theoretical Probability} 23(4), 1227-1243.


\bibitem{JiangQi}
Jiang, T. and Qi, Y. (2017). Spectral radii of large non-Hermitian random matrices. {\it J.  Theor. Probab.} 30, 326-364.

%\bibitem{Jo}
%Johansson, K. (2000). Shape fluctuations and random matrices. {\it Comm. Math. Phys.} 209, 437-476.

\bibitem{Johansson07}
Johansson, K. (2007). From Gumbel to Tracy-Widom. {\it Probab. Theory Relat. Fields} 138, 75-112.

\bibitem{John}
Johnstone, I. (2001). On the distribution of the largest eigenvalue in principal compo-
nents analysis. {\it Ann. Stat.} 29, 295-327.

\bibitem{Johnstone}
Johnstone, I. (2008). Multivariate analysis and Jacobi ensembles: Largest eigenvalue, Tracy–Widom limits and rates of convergence. {\it Ann. Stat.}, 36(6), 2638-2716.


%\bibitem{Khor}
%Khoruzhenko, B. A. and  Sommers, H. J. (2001). Non-Hermitian Random Matrix Ensembles. In ``The Oxford Handbook of Random Matrix Theory (Oxford %Handbooks in Mathematics)". {\it Oxford University Press}, 376-397.

%\bibitem{Kostlan}
%Kostlan, E. (1992). On the spectra of Gaussian matrices. {\it Linear Algebra and Its Applications }162-164, 385-388.

%\bibitem{Kris}
%Krishnapur, M. (2009). From random matrices to random analytic functions. {\it Ann. Probab. 37(1), 314-346}.


%\bibitem{Ku}
%Kuijlaars, A. B. J.  and L\'{o}pez-Garc\'{i}a, A. (2015).  The normal matrix model with a monomial potential, a vector equilibrium problem, %and multiple orthogonal polynomials on a star. {\it Nonlinearity} 28, 347-406.
%
%A vector equilibrium problem for the normal matrix model,
%and multiple orthogonal polynomials on a star. Available at {\tt http://arxiv.org/pdf/1401.2419v1.pdf}.


%\bibitem{Lehman}
%Lehmann, N. and Sommers, H. J. (1991). Eigenvalue statistics of random real matrices. {\it Phys.
%Rev. Lett.} 67, 941-944.


\bibitem{Mehta}
Mehta, M. L. (2004). {\it Random matrices}. Volume 142. Academic Press.

\bibitem{MS} Mezzadri, F. and Snaith, N. C. (2005). {\it Recent perspectives in random matrix theory and number theory}.
Cambridge Univ Press.


%\bibitem{Rourke}
%O'Rourke, S. and Soshnikov, A.  (2011). Products of independent non-Hermitian random matrices. {\it Electrical Journal of Probability} 16(81), %2219-2245.

%\bibitem{Rourke14}
%O'Rourke, S., Renfrew, D., Soshnikov, A. and Vu, V.  (2014). Products of independent elliptic random matrices. Available at {\tt %http://arxiv.org/pdf/1403.6080v2.pdf}.

%\bibitem{Petrov}
%Petrov, V.V. (1975). Sums of Independent Random Variables. {\it Springer-Verlag}.


%\bibitem{Petz}
%Petz, D. and Hiai, F. (1998). Logarithmic energy as entropy functional.  In Advances in Differential Equations and Mathematical Physics, eds.  %Carlen, Harrell, Loss. {\it Contemporary Math.} 217, 205-221.

%\bibitem{Pillai}
%Pillai, N. S. and Yin, J. (2014). Universality of covariance matrices. {\it Ann. Appl. Probab.} 24(3), 935-1001.

%\bibitem{Ramirez}
%Ram\'{\i}rez, J., Rider, B. and Vir\'{a}g, B. (2011). Beta ensembles, stochastic Airy
%spectrum and a diffusion. {\it J. Amer. Math . Soc.} 24, 919-944.

\bibitem{Reiss}
Reiss, R. D. (1981).  Uniform approximation to distributions of extreme order statistics.
{\it Advances in Applied Probability} 13, 533-547.


%\bibitem{Resnick}
%Resnick, S. I. (2007). Extreme Values, Regular Variation and Point Processes. {\it Springer}.


\bibitem{Rider}
Rider, B. C. (2003). A limit theorem at the edge of a non-Hermitian random matrix ensemble. {\it J. Phys. A}  36(12), 3401-3409.

\bibitem{Rider09_23}
Rider, B. C. (2004). Order statistics and Ginibre’s ensembles. {\it Journal of Statistical Physics} 114, 1139-1148.


\bibitem{RS}
Rider, B. C. and Sinclair, C. D. (2014). Extremal laws for the real Ginibre ensemble. {\it Ann. Appl. Probab.} 24(4), 1621-1651.

\bibitem{TaoVu}
Tao, T. and Vu, V. (2011). Random matrices: Universality of local eigenvalue statistics. {\it Acta Mathematica} 206(1), 127-204.

%\bibitem{Tao}
%Tao, T. and Vu, V. (2014). Random matrices: The universality phenomenon for Wigner ensembles. In Modern Aspects of Random Matrix Theory %(Proceedings of Symposia in Applied Mathematics, ed. Vu, V.). {\it American Mathematical Society}.

\bibitem{Tracy94}
Tracy, C. A. and Widom, H. (1994). Level-spacing distributions and Airy kernal. {\it Comm.
Math. Physics} 159, 151-174.

\bibitem{Tracy96}
Tracy, C. A. and Widom, H. (1996). On the orthogonal and symplectic matrix ensembles.
{\it Comm. Math. Physics} 177, 727-754.

\bibitem{TW02}
Tracy, C. A. and Widom, H. (2002).  Distribution functions for largest eigenvalues and their applications. {\it Proceedings of the
ICM, Beijing} 1, 587-596.


%\bibitem{Watcher}
%Wachter, K. W. (1980). The limiting empirical measure of multiple discriminant ratios. {\it Ann.
%Statist.} 8, 937-957.


\bibitem{Wigner}
Wigner, E. P. (1955). Characteristic vectors of bordered matrices with infinite dimensions. {\it Ann. Math.} 62, 548-564.

\bibitem{Wishart}
Wishart, J. (1928). The generalized product moment distribution in samples
from a normal multivariate population. {\it Biometrika} 20, 35-52.

\bibitem{Zski}
\.{Z}yczkowski, K. and Sommers, H. (2000). Truncation of random unitary matrices. {\it J. Phys. A: Math. Gen.} 33, 2045-2057.
\end{thebibliography}
\end{document}